\documentclass[10pt,reqno]{amsart}
\usepackage{latexsym}
\usepackage{graphicx}
\usepackage{fancyhdr,amssymb}

\usepackage{hyperref} 

\title{Sequences, Series and Uniform distribution of Square-Prime numbers}
\author{Raghavendra N Bhat}
\address{\newline Raghavendra N Bhat\newline 
University of Illinois, Urbana Champaign \newline
Department of Mathematics\newline
1409 West Green Street\newline
Urbana, IL 61801
}
\email{rnbhat2@illinois.edu (Corresponding Author)}
\begin{document}
\maketitle
\begin{abstract}


We defined numbers of the form $p\cdot a^2$ as \textbf{SP numbers} (Square-Prime numbers) ($a\neq1$, $p$ prime) in \cite{Bhat} along with proofs on their distribution. These numbers are listed in the OEIS as A228056 \cite{OEIS}. Some examples of SP numbers : 75 = 3 $\cdot$ 25; 108 = 3 $\cdot$ 36; 45 = 5 $\cdot$ 9. This paper explores sequences of these numbers, sequences between 0 and 1 related to these numbers and analyzes the distribution of some of these sequences.
\end{abstract}
\section{Introduction and Summary of previous results}
In our previous paper titled `Distribution of square-prime numbers', we introduced the notion of a product of a prime and a square, with the square not being 1. We call these SP numbers.\\\par Here are the first 25 SP numbers:
8, 12, 18, 20, 27, 28, 32, 44, 45, 48, 50, 52, 63, 68, 72, 75, 76, 80, 92, 98, 99, 108, 112, 116, 117.\\\par In the paper we mainly proved three results. Firstly, SP numbers have an asymptotic distribution similar to prime numbers, meaning, for a large natural number $n$, the number of SP numbers smaller than $n$ is asymptotic to $\frac{n}{\log n}$, as with the case of primes. We denote the number of SP numbers smaller than $n$ by $SP(n)$. The result derived in the paper is $SP(n) = (\zeta(2) - 1)\cdot \frac{n}{\log n} + O\left(\frac{n}{\log n^2}\right).$ We see that although asymptotically similar to the distribution of primes, the actual value of $SP(n)$ is lower than $\pi(n)$ for large $n$, since $\zeta(2) - 1 < 1$.
\par Secondly, there are infinitely many pairs of SP numbers that are consecutive natural numbers. We call these SP twins. An example is (27,28).
\par Thirdly, we analyzed techniques to give asymptotic estimates for number of SP numbers having a certain last digit. For example, the number of SP numbers ending in 1 approaches $\frac{1}{400}\frac{n}{\log n}(\zeta(2,1/10) + \zeta(2,9/10) + \zeta(2,3/10) + \zeta(2,7/10) - 4).$ \par In this paper, our goal is to study sequences of SP numbers, and how these sequences are distributed. Sometimes, we study the distribution (uniform) by restricting our numbers between 0 and 1.
\section{Equidistribution of SP numbers}
Let $U = \{ U_1, U_2, U_3, U_4,...\}$ be an infinite sequence of finite sequences with each individual sequence containing numbers between 0 and 1. We define $U$ to be \textbf{uniformly distributed} if and only if for any $0\leq \alpha \leq \beta \leq 1,$
\[ \lim_{j\to\infty} \#\{k \in (\alpha, \beta) | k \in U_j\} = \beta - \alpha, \] where $\#$ denotes the number of elements in the set.\\\par Now consider $U=\{ U_1, U_2, U_3, U_4, \cdot\}$, where each $U_j=\{\frac{sp_1}{j}$, $\frac{sp_2}{j},\dots,\frac{sp_{n_j}}{j}\}$, with $sp_1$, $sp_2$, \dots $sp_{n_j}$ being SP numbers in order starting from 8 and ending at the largest SP lesser than $j$. We represent the largest SP number lesser than $j$ by $n_j$. Thus, in every $U_j$, all numbers are restricted to the interval $(0,1)$.\\\\
\textbf{Theorem 2.1 :} $U$ is uniformly distributed.
\begin{proof}
Fix $0 \leq \alpha < \beta \leq 1$.\\\\
We explore the following limit:
\[ \lim_{j\to\infty} \frac{\# \{ 1 \leq k \leq n_j | \frac{sp_k}{j} \in [\alpha, \beta]\}}{SP(j)},\] essentially looking at the fraction of numbers that are between $\alpha$ and $\beta$ as compared to the total numbers in $U_j$, which is the number of SP numbers smaller than $j$, i.e. $SP(j)$. The numerator of the limit
\[=\# \{ 1 \leq k \leq n_j | \frac{sp_k}{j} \in [\alpha, \beta]\} = \# \{ 1 \leq k \leq n_j | sp_k \in [\alpha j, \beta j] \Rightarrow SP(\beta j) - SP(\alpha j)\}.\] Thus, our intended limit can be calculated as follows:
$$\lim_{j\to\infty} \frac{SP(\beta j) - SP(\alpha j)}{SP(j)} = \lim_{j\to\infty} \frac{(\zeta - 1)\cdot \frac{\beta j}{\log \beta j} + O\left(\frac{\beta j}{\log^2 \beta j}\right) - (\zeta - 1)\cdot \frac{\alpha j}{\log \alpha j} + O\left(\frac{\alpha j}{\log^2 \alpha j}\right)}{(\zeta - 1)\cdot \frac{j}{\log j} + O\left(\frac{j}{\log^2 j}\right)}$$
\\
$$= \lim_{j\to\infty} \frac{(\zeta - 1)\cdot \frac{\beta j}{\log \beta j}+ O\left(\frac{\beta j}{\log^2 \beta j}\right) - \frac{\alpha j}{\log \alpha j} + O\left(\frac{\alpha j}{\log^2 \alpha j}\right)}{(\zeta - 1)\cdot \frac{j}{\log j} + O\left(\frac{j}{\log^2 j}\right)}$$
\\
$$= \lim_{j\to\infty} \frac{\frac{\beta j}{\log j}\cdot \left(1 + O\left(\frac{1}{\log j}\right)\right) - \frac{\alpha j}{\log j}\cdot \left(1 + O\left(\frac{1}{\log j}\right)\right)}{\frac{j}{\log j}\cdot (1 + O\left(\frac{1}{\log j}\right)}$$
\\
$$= \lim_{j\to\infty} \frac{\beta\cdot \left(1 + O\left(\frac{1}{\log j}\right)\right) - \alpha\cdot \left(1 + O\left(\frac{1}{\log j}\right)\right)}{(1 + O\left(\frac{1}{\log j}\right)} = \beta - \alpha$$
\\
Thus, $U$ is uniformly distributed by definition.
\end{proof}
\section{Farey-like SP sequence}
We now study a sequence based on ideas from the well-known Farey sequence in number theory. The Farey sequence $F_n$ for any positive integer $n$ is the set of irreducible rational numbers $\frac{a}{b}$, arranged in increasing order where $0\leq a\leq b\leq n$, with $a,b$ being co-prime, \cite{farey}.\\\par Consider large $x$. Create a sequence SPFarey$_x$ where $\frac{sp_1}{sp_2} \in$ SPFarey$_x$ if $sp_1 < sp_2$, $(sp_1,sp_2)=1$ and $sp_1, sp_2$ are both SP numbers smaller than or equal to $x$. SPFarey$_x$ is thus a sequence of ratios of SP number pairs upto $x$ that are co-prime. As an example, SPFarey$_{50}$ = $\{\frac{8}{27}, \frac{8}{45}, \frac{20}{27}, \frac{27}{28}, \frac{27}{32}, \frac{27}{44}, \frac{27}{50}, \frac{28}{45}, \frac{32}{45}, \frac{44}{45}\}.$\\\\
\textbf{Theorem 3.1:} Define SPFarey($x$) as the cardinality of SPFarey$_x$. Then, SPFarey($x$) = $O\left(\frac{x^2}{\log x^2}\right)$ for large $x$. \begin{proof} We have an asymptotic formula \cite{Bhat} for $SP(x)$ i.e. $(\zeta(2) - 1)\cdot\frac{x}{\log x} + O\left(\frac{x}{\log^2 x}\right)$. Thus, the total number of pairs of SP numbers smaller than $x$ is
$$(\zeta(2) - 1)^2\cdot\frac{x^2}{\log^2 x} + O\left(\frac{x^2}{\log^4 x}\right) + 2\left((\zeta(2) - 1)\cdot\frac{x}{\log x}\right) + O\left(\frac{x}{\log^2 x}\right)$$
We divide the equation by 2 to avoid duplicate counting of ratios (our sequence restricts the numerator to be smaller than the denominator).
$$= \frac{(\zeta(2) - 1)^2}{2}\cdot\frac{x^2}{\log^2 x} + \left((\zeta(2) - 1)\cdot\frac{x}{\log x}\right) + O\left(\frac{x}{\log^2 x}\right)$$
The probability that a randomly chosen pair of numbers is co-prime is $\frac{6}{\pi^2} = \frac{1}{\zeta(2)}$. We proved in Theorem 2.1 that SP numbers are uniformly distributed, hence allowing us to calculate an asymptotic form of SPFarey$(x)$
$$\frac{(\zeta(2) - 1)^2}{2\cdot\zeta(2)}\cdot\frac{x^2}{\log^2 x} + \left(\frac{\zeta(2) - 1}{\zeta(2)}\cdot\frac{x}{\log x}\right) + O\left(\frac{x}{\log^2 x}\right).$$
Setting $A=(\zeta(2) - 1)^2$ we have our asymptotic estimate for SPFarey$(x):$
$$\frac{A}{2\cdot\zeta(2)}\cdot\frac{x^2}{\log^2 x} + \left(\frac{A^{\frac{1}{2}}}{\zeta(2)}\cdot\frac{x}{\log x}\right) + O\left(\frac{x}{\log^2 x}\right)$$
$$ = \frac{1}{\zeta(2)}\left(\frac{Ax^2}{2\log^2 x} + \frac{A^{\frac{1}{2}}x}{\log x} + \right) + O\left(\frac{x}{\log^2 x}\right).$$
\end{proof}
\section{SP Harmonic Series}
The harmonic series of natural numbers is defined as the following infinite sum:
$$\sum_{n=1}^{\infty} \frac{1}{n}$$
It is well known that for a fixed number $k$, the partial sum of the harmonic series, up to $k$ is asymptotic to $\log k$
$$\sum_{n=1}^{k} \frac{1}{n} = \log k + O(1)$$
We will now explore a variant of this series for SP numbers. Consider the SP Harmonic series defined as
$$\sum_{i=1}^{\infty} \frac{1}{sp_i}$$ where $sp_i$ is the $i^{th}$ SP number. Here are the first few terms of the summation:
$$\frac{1}{8} + \frac{1}{12} + \frac{1}{18} + \frac{1}{20} + \frac{1}{27} + \frac{1}{28} + \frac{1}{32} + ....$$
\textbf{Theorem 2.1 :} The SP Harmonic series diverges.
\begin{proof}
We will soon explore an asymptotic estimate for the partial sums of the SP Harmonic Series. However, to prove that it diverges is straightforward and does not require any analytic tools.\\\\
Consider the prime harmonic series:
$$\sum_{i=1}^{\infty}\frac{1}{p_i}$$ where $p_i$ is the $i^{th}$ prime number. We know from Euler \cite{Mathworld} that the prime harmonic series diverges.\\\\
Consider
$$S = \frac{1}{4} \sum_{i=1}^{\infty}\frac{1}{p_i} = \frac{1}{4}\left(\frac{1}{2} + \frac{1}{3} + \frac{1}{5} + ....\right)$$
$$= \frac{1}{8} + \frac{1}{12} + \frac{1}{20} + ...$$
We know $S$ diverges as it is a constant multiple of the prime harmonic series. Closer inspection tells us that $S \leq$ SP Harmonic series because it is only considering SP numbers of the form $4p^2$. Thus, the SP Harmonic diverges.\end{proof}
\noindent\textbf{Theorem 2.2 :} Define SPHarmonic($k$) as
$$\sum_{i=1}^{t} \frac{1}{sp_i}$$ where $sp_t$ is the largest SP number smaller than or equal to $k$. Then SPHarmonic($k$) $= O\left(\log (\log k) \right)$
\begin{proof}
Since SPHarmonic($k$) is a sum over all reciprocals of SP numbers smaller than or equal to $k$, we represent it as the following double sum:
$$\sum_{1<a\leq\sqrt{\frac{k}{2}}}\sum_{p\leq \frac{k}{a^2}} \frac{1}{p}$$
Simplifying based on Meissel-Merten and Euler, we get
\begin{equation}\label{1}
\sum_{1<a\leq\sqrt{\frac{k}{2}}} \frac{1}{a^2}\left(\log \log \left(\frac{k}{a^2}\right) + M\right)
\end{equation}with $M$ being the Meissel-Merten constant \cite{MM}. We can break up (1) into two summations $S_1$ and $S_2$ as follows:
$$S_1 = \sum_{1<a\leq A} \frac{1}{a^2}\left(\log \log \left(\frac{k}{a^2}\right) + M\right)$$
$$S_2 = \sum_{A\leq a\leq \sqrt{\frac{k}{2}}} \frac{1}{a^2}\left(\log \log \left(\frac{k}{a^2}\right) + M\right)$$ where $A$ is fixed. We first simplify $S_2$, the tail:
$$S_2 = O\left(\sum_{A\leq a\leq \sqrt{\frac{k}{2}}} \frac{1}{a^2}\left(\log \log \left(k\right) + M\right)\right) < O\left(\sum_{A\leq a\leq \infty} \frac{1}{a^2}\left(\log \log \left(k\right) + M\right)\right)$$
$$=(\log \log k + M)\cdot O\left(\int_{A}^{\infty} \frac{1}{a^2} \,da\right)=(\log \log k + M)\cdot O\left(\frac{1}{A}\right)=O\left(\frac{\log \log k}{A}\right)$$\\
Now, simplifying $S_1$, the main term:
$$\sum_{1<a\leq A} \frac{1}{a^2}\left(\log \log \left(\frac{k}{a^2}\right) + M\right)$$
$$\log\left(\log\left(\frac{k}{A^2}\right)\right) = \log\left(\log k - 2\log A\right) = \log\log k + \log\left(1 - \frac{2\log A}{\log k}\right)$$\\
Using the Taylor series expansion of $\log(1-x)$ for small $x$, and also since $A\geq a$, we get
$$\log\log\left(\frac{k}{a^2}\right) = \log\log k + O\left(\frac{2\log A}{\log k}\right)$$
$$\Rightarrow S_1 = \sum_{1<a\leq A} \frac{1}{a^2}\left(\log\log k + O\left(\frac{2\log A}{\log k}\right) + M\right)$$
$$= \log \log k\sum_{1<a\leq A}\frac{1}{a^2} + O\left(\sum_{1<a\leq A}\frac{1}{a^2}\cdot\frac{2\log A}{\log k}\right) + \sum_{1<a\leq A}\frac{1}{a^2}M$$\\
Combining $S_1$ and $S_2$, we get SPHarmonic($k$) as $k \to \infty$
$$S_1 + S_2 = \log \log k\sum_{1<a\leq A}\frac{1}{a^2} + \sum_{1<a\leq A}\frac{1}{a^2}M + O\left(\frac{\log A}{\log k}\right) + O\left(\frac{\log \log k}{A}\right)$$
Let $A = (\log k)^\beta$, with $\beta$ an arbitrary natural number.
\begin{equation}
\text{SPHarmonic}(k) = \left(\frac{\pi^2}{6} - 1\right)(\log \log k + M) + O\left(\frac{\log (\log k)^\beta}{\log k}\right) + O\left(\frac{\log \log k}{(\log k)^\beta}\right)
\end{equation}
\end{proof}
\section{Comparison Table of SPHarmonic(x)}
Here is a table comparing the actual value of SPHarmonic($x$) with the dominant term in the asymptotic estimate from (2). One immediately notices that the error terms are negligible as for large $X$.
\begin{center}
\begin{tabular}{||c|c|c||} 
 \hline
 X & SPHarmonic(X) & $\left(\frac{\pi^2}{6} - 1\right)(\log\log X)$\\ [0.5ex] 
 \hline\hline
 100 & 0.6375 & 0.9849 \\ 
 \hline
 1000 & 1.0355 & 1.2464 \\
 \hline
 10000 & 1.3237 & 1.4319 \\
 \hline
 100000 & 1.5342 & 1.5759\\
 \hline
 250000 & 1.6029 & 1.6253\\
 \hline
\end{tabular}
\end{center}
\section{Future work and Conjectures}
One wonders if the partial sums of the harmonic series of SP twins would converge and, if so, have any relationship to Brun's constant \cite{brun}. Formally, the problem of interest is the asymptotic estimate of the following summation:
$$\sum_{s,s+1 \in SP}\frac{1}{s} + \frac{1}{s+1}$$
An observation that can be made based on the way the summation is set up is the allowance of duplicate counting of SP reciprocals if it shows up in two twin SP pairs. Thus, if $x,y,z$ are all SP numbers and consecutive, our summation would involve $\frac{1}{x} + \frac{1}{y} + \frac{1}{y} + \frac{1}{z}$, owing to the fact that $y$ is part of two paris of SP twins.
\section{Acknowledgements} I would like to acknowledge the professors at the University of Illinois for their valuable advice and support. In particular, I am grateful to Professors Alexandru Zaharescu and Bruce Berndt for their guidance and inputs.
\section{Data Availability Statement}
Data sharing not applicable to this article as no data sets were generated or analysed during the current study.
\section{Conflict of Interest}
The author declares that there is no conflict of interest.

\end{document}